\documentclass[a4paper,12pt]{article}
\usepackage{amssymb,latexsym,amsmath,amsfonts}
\usepackage{mathrsfs}
\usepackage{amsthm}
\usepackage{bm}
\usepackage[mathcal]{euscript}
\usepackage{bbm}

\newcommand{\cF}{\mathscr{F}}
\newcommand{\R}{\mathbb{R}}

 \newcommand{\EE}{{\sf E}}
 
 \newcommand{\PP}{{\sf P}}

\renewcommand{\le}{\leqslant}
\renewcommand{\ge}{\geqslant}

\newcommand{\ind}{{\mathbbm 1}}

 \newtheorem{thm}{\indent Theorem}
 \newtheorem{cor}{\indent Corollary}

\begin{document}

                    \begin{center}
                    {\Large\textbf{On pathwise counterparts of Doob's maximal inequalities{\normalfont\normalsize\footnote{This research was supported by the Russian Scientific Fund (project N 14-21-00162)}}}}
                    \end{center}

                    \begin{center}
                    {\Large \textbf{A.\,A.\,Gushchin{\normalfont\normalsize\footnote{Steklov Mathematical Institute, Moscow, Russia. E-mail: gushchin@mi.ras.ru}} } }
                    \end{center}
 \begin{abstract}
 {\small In this paper, we present pathwise counterparts of Doob's maximal inequalities (on the probability of exceeding a level) for submartingales and supermartingales.
}
\end{abstract}

Recently a new method of proving martingale inequalities became popular. Namely, they are derived (and thus often improved) from elementary \emph{deterministic} inequalities. As typical examples, let us mention Doob's maximal $L^p$- and $L\log L$-inequalities \cite{ABPST13} and the Burkholder--Davis--Gundy inequality \cite{BeiSio13}. Deterministic inequalities have a natural interpretation in terms of robust hedging of options, see \cite{Hob11}, which served as the impetus for their appearance. We also mention the papers \cite{OblYor06}, \cite{BouNutz13}, \cite{{BeiNutz14}}, \cite{CoxObl14}, \cite{H-LabOblSpoTou14}, \cite{OblSpoTou14}, dealing with similar problems.

The purpose of this note is to present elementary pathwise counterparts of Doob's maximal inequalities on the probability of exceeding a level. Substituting a trajectory of a stochastic process in our inequality and taking expectations, we obtain Doob's inequality for supermartingales and submartingales, see \cite[Chapter~VII, Theorem~3.2]{doob53}, due to the fact that one of the terms in the inequalities can be dropped if the process is a supermartingale (in the case of the first inequality) or a submartingale (in the case of the second one). We also show that the pathwise counterparts of Doob's maximal $L^p$- and $L\log L$-inequalities from the paper \cite{ABPST13} can be obtained by integration from our inequalities.

Let $x=(x_0,\dots,x_n)$ be a vector of real numbers. Put
\[
\bar x_k = \max\,\{x_0,\dots,x_k\}, \quad k=0,\dots n;\qquad
\Delta x_k = x_k-x_{k-1}, \quad k=1,\dots n.
\]
The symbol $\ind_A$ stands for the indicator function that is $1$ on a set $A$ and $0$ outside $A$.

\begin{thm}
For any $\lambda\in\R,$
\begin{align}
\lambda\ind_{\{\bar x_n\ge\lambda\}} & \le x_0\wedge\lambda +\sum_{k=1}^n \ind_{\{\bar x_{k-1}<\lambda\}}\Delta x_k - x_n\ind_{\{\bar x_n<\lambda\}},\label{e:1}\\
\lambda\ind_{\{\bar x_n\ge\lambda\}} &\le -(x_0-\lambda)\ind_{\{x_0\ge\lambda\}} -\sum_{k=1}^n \ind_{\{\bar x_{k-1}\ge\lambda\}}\Delta x_k  + x_n\ind_{\{\bar x_n\ge\lambda\}}\label{e:2}.
\end{align}
\end{thm}

The proof of inequalities \eqref{e:1} and \eqref{e:2} reduces to their trivial verification in three cases: $\bar x_n<\lambda$, $x_0\ge\lambda$ è $\bar x_{j-1}<\lambda\le x_j$ ($j=1,\dots n$). Moreover, in the first two cases both inequalities are equalities and, in the third case, the difference of the right-hand and left-hand sides equals $x_j-\lambda$ in both inequalities.

Now let $X=(X_k)_{k=0,1,\dots,n}$ be a stochastic process given on a filtered probability space $(\Omega,\cF,(\cF_k)_{k=0,1,\dots,n},\PP)$.

\begin{cor} Let $\lambda\in\R$.

{\rm (i)} If $X$ is a supermartingale{\rm ,} then
\begin{equation}\label{e:3}
\lambda\PP(\bar X_n\ge\lambda) \le \EE(X_0\wedge\lambda) - \int_{\{\bar X_n<\lambda\}}X_n\,d\PP.
\end{equation}

{\rm (ii)} If $X$ is a submartingale{\rm ,} then
\begin{equation}\label{e:4}
\lambda\PP(\bar X_n\ge\lambda) \le -\EE\bigl[(X_0-\lambda)\ind_{\{X_0\ge\lambda\}}\bigr] + \int_{\{\bar X_n\ge\lambda\}}X_n\,d\PP.
\end{equation}
\end{cor}

Inequalities \eqref{e:3} and \eqref{e:4} slightly improve the original inequalities due to Doob  \cite[Chapter~VII, Theorem~3.2]{doob53} that are obtained if we replace $\EE(X_0\wedge\lambda)$ by $\EE(X_0)$ in \eqref{e:3} and drop the first term $-\EE\bigl[(X_0-\lambda)\ind_{\{X_0\ge\lambda\}}\bigr]$ on the right in \eqref{e:4}.

Now let all $x_0,\dots,x_n$ be nonnegative, $p>1$, $q=p/(p-1)$. Then the following relations hold, where the first inequality  follows from \eqref{e:2} and the second one follows from the inequality  $ab\le a^p/p + b^q/q$ ($a,b\ge 0$):
\begin{align*}
\bar x_n^p & = p\int\limits_0^\infty \lambda^{p-1}\ind_{\{\bar x_n\ge\lambda\}}\,d\lambda\\
           &\le p\int\limits_0^\infty \lambda^{p-2}x_n\ind_{\{\bar x_n\ge\lambda\}} \,d\lambda -p\int\limits_0^\infty \lambda^{p-2}(x_0-\lambda)\ind_{\{x_0\ge\lambda\}}\,d\lambda - p\sum_{k=1}^n \int\limits_0^\infty \lambda^{p-2} \ind_{\{\bar x_{k-1}\ge\lambda\}}\Delta x_k \,d\lambda\\
           &= qx_n\bar x_n^{p-1} - qx_0^p + x_0^p -q\sum_{k=1}^n \bar x_{k-1}^{p-1}\Delta x_k\\
           &\le \frac{\bar x_n^{p}}{q} + \frac{q^p x_n^p}{p} - (q-1)x_0^p -q\sum_{k=1}^n \bar x_{k-1}^{p-1}\Delta x_k.
\end{align*}
Therefore,
\begin{equation}\label{e:5}
\bar x_n^p \le q^p x_n^p - qx_0^p -qp\sum_{k=1}^n \bar x_{k-1}^{p-1}\Delta x_k.
\end{equation}
Inequality \eqref{e:5} was obtained in \cite{ABPST13}. It implies the following minor generalization of Doob's maximal $L^p$-inequality: if $X$ is a nonnegative submartingale and $\EE X_n^p<\infty$, then
\[
\EE\bigl[\bar X_n^p\bigr]  \le q^p\EE [X_n^p]-q\EE [X_0^p],
\]
see \cite{Cox84}.

In conclusion let us consider the case $p=1$. Assume additionally that $x_0>0$. Then the following relations hold, where the first inequality  follows from \eqref{e:2} and the second one follows from the inequality $a\log\,b\le a\log\,a + e^{-1}b$ ($a\ge 0$, $b>0$):
\begin{align*}
\bar x_n & = x_0 + \int\limits_{x_0}^\infty \ind_{\{\bar x_n\ge\lambda\}}\,d\lambda\\
           &\le x_0 + x_n\int\limits_{x_0}^\infty \lambda^{-1}\ind_{\{\bar x_n\ge\lambda\}} \,d\lambda - \sum_{k=1}^n \Delta x_k \int\limits_{x_0}^\infty \lambda^{-1} \ind_{\{\bar x_{k-1}\ge\lambda\}}\,d\lambda\\
           &= x_0 + x_n\log \bar x_n - x_n\log x_0 - \sum_{k=1}^n \log (\bar x_{k-1}/x_0) \Delta x_k \\
           &\le x_0+x_n\log (x_n/x_0)+e^{-1}\bar x_n - \sum_{k=1}^n \log (\bar x_{k-1}/x_0)\Delta x_k.
\end{align*}
Therefore,
\begin{align}
\bar x_n &\le \frac{e}{e-1}\Bigl[x_0+ x_n\log (x_n/x_0) - \sum_{k=1}^n \log (\bar x_{k-1}/x_0)\Delta x_k\Bigr]\label{e:6}\\
&= \frac{e}{e-1}\Bigl[x_0(1-\log x_0) + x_n\log x_n - \sum_{k=1}^n \log\, \bar x_{k-1}\Delta x_k\Bigr].\label{e:7}
\end{align}
This inequality with the right-hand side as in \eqref{e:7} is proved in \cite {ABPST13}. Rewriting it in the form \eqref{e:6} allows us to drop the last term with the sum after substituting a nonnegative submartingale $X$ for $x$ and taking expectations. At the same time, in general, the last term in \eqref{e:7} has an indefinite sign after substituting a submartingale and taking expectations if $X$ is not a martingale. More precisely, if $X$ is a nonnegative martingale such that $\EE [X_n\log\,X_n]<+\infty$, then the inequality
\begin{equation}\label{e:8}
\EE[\bar X_n] \le \frac{e}{e-1}\bigl\{\EE[X_0(1-\log X_0)] + \EE[X_n\log\, X_n]\bigr\}
\end{equation}
holds. It is easy to see that this is not true in general if $X$ is a (strictly positive) submartingale: it is enough to put $n=1$, $X_0=\varepsilon$, where $\varepsilon>0$ is small enough, and $X_1=1$. In other words, inequality (Doob-$L^1$) in the statement of Theorem~1.1 in~\cite{ABPST13} is valid for nonnegative martingales and is not valid for submartingales as is stated in this theorem. Nevertheless, the following improvement of Doob's maximal $L\log\,L$-inequality is true: for any nonnegative submartingale $X$,
\begin{equation}\label{e:9}
\EE[\bar X_n] \le \frac{e}{e-1}\bigl\{1 + \EE[X_n\log\, X_n]\bigr\}.
\end{equation}
For martingales, \eqref{e:9} follows from \eqref{e:8}. If $X$ is a submartingale, then the inequality follows from \eqref{e:9} applied to the martingale $Y_k=\EE[X_n|\cF_k]$, since, clearly, $Y_n=X_n$ and $\bar Y_n \ge \bar X_n$. Recall that, in Doob's maximal $L\log\,L$-inequality, unlike \eqref{e:9}, there appears $\log^+$ instead of $\log$, and that the constant $e/(e-1)$ in \eqref{e:9} is optimal, see \cite{Gil86}.

\end{document}